\documentclass[11pt, a4paper, twoside]{amsart}

\usepackage{amsmath}
\usepackage{amsfonts} 
\usepackage{amssymb}
\usepackage{latexsym}

\usepackage{hyperref}
\hypersetup{backref, pdfpagemode=FullScreen, colorlinks=true,
  citecolor=magenta, linkcolor=cyan, urlcolor=blue}

\usepackage{mathtools}
\mathtoolsset{showonlyrefs}

\usepackage{graphicx}
\usepackage{xcolor}
\usepackage{graphicx}

\theoremstyle{change}
\newtheorem{theorem}{Theorem}[section]
\newtheorem*{theo*}{Theorem}

\newtheorem{corollary}[theorem]{Corollary}

\newtheorem{conjecture}[theorem]{Conjecture}

\newtheorem{definition}[theorem]{Definition}

\newcommand{\bad}{B^d}
\newcommand{\sad}{\mathbb{S}^{d-1}}
\newcommand{\Kd}{\mathcal{K}^d}
\newcommand{\packk}{{\mathcal P}(K)}
\newcommand{\packkn}{{\mathcal P}_n(K)}

\newcommand{\packbn}{{\mathcal P}_n(B^d)}
\newcommand{\deninf}{\delta(K)}
\newcommand{\vol}{\mathrm{vol}\,}
\newcommand{\R}{\mathbb{R}} 
\newcommand{\N}{\mathbb{N}} 
\newcommand{\Z}{\mathbb{Z}}
\newcommand{\ip}[2]{\left\langle #1,#2\right\rangle}

\newcommand{\vx}{{\boldsymbol x}}
\newcommand{\vy}{{\boldsymbol y}}
\newcommand{\vu}{{\boldsymbol u}}
\newcommand{\vt}{{\boldsymbol t}}
\newcommand{\V}{\mathrm{V}}
\newcommand{\vnull}{{\boldsymbol 0}}
\newcommand{\bd}{\mathrm{bd}\,}
\newcommand{\inte}{\mathrm{int}\,}
\newcommand{\aff}{\mathrm{aff}\,}
\newcommand{\conv}{\mathrm{conv}\,}
\newcommand{\ov}{\overline}
\newcommand{\trans}{\intercal}

\newcommand{\norm}[2]{\left\vert #1\right\vert_{#2}}
\newcommand{\enorm}[1]{\left\lvert #1\right\rvert}

\numberwithin{equation}{section}

\begin{document}

\title{Packings, sausages and catastrophes}
\author{Martin Henk}
\address{Technische Universität Berlin, Institut für Mathematik, Sekr. MA4-1, Stra{\ss}e des 17 Juni 136, D-10623 Berlin}
\email{henk@math.tu-berlin.de}
\author{Jörg M. Wills}
\address{Universität Siegen, Emmy-Noether-Campus, Walter-Flex-Str. 3, 
  D-57068 Siegen}
\email{wills@mathematik.uni-siegen.de}
\dedicatory{Dedicated to the memory of Uli Betke (1948--2008)}
\begin{abstract} In this survey we give an  overview about
  some of the main results on parametric densities,
  a concept which unifies the theory of finite (free) packings
  and the classical theory of infinite packings.      
\end{abstract}
\maketitle

\section{Introduction}
The theory of infinite packings of convex bodies, in particular, lattice packings
of spheres is a fundamental and classical topic in mathematics which plays a role in various
branches of mathematics as  number theory, group theory, geometry of numbers, algebra,
and which has numerous  applications to coding theory, cryptography, crystallography and more.
Here the main problem is to arrange infinitely many  non-overlapping (lattice-) translative copies of 
a given convex body such that the whole space is covered as much as
possible, i.e.,  packed as densely a possible. 

On the other hand one may say, that all packings in real world are finite,
even the atoms in crystals or sand at the beach, and in the theory of finite packings we want to arrange
 finitely many  non-overlapping (lattice-) translative copies of a convex body as good as possible. 
Roughly speaking, there are two different ways to specify 
``as good as possible'':    
The first one leads to so called {\em bin packings} where a container (bin) of a prescribed
shape (ball, simplex, cube, etc.) but of minimal size (volume) is looked for containing a given number of non-overlapping (lattice-) translative copies of the given convex body. In contrast to this,
we consider here so called {\em free (finite) packings} where the goal is to minimize the volume of the convex hull of a given number of non-overlapping (lattice-) translative copies of the convex body.    
 
The volume--based free packing approach was introduced in 1892 by the Norwegian
mathematician Axel Thue. 
 For a given number of circles, he
considered all possible packings and their convex hulls and
asked for the minimal volume of these convex hulls. 
Thue's approach was further developed between
1940 and 1972 by many prominent
mathematicians, e.g., by ~L.~Fejes T\'oth, R.P.~Bambah, 
 C.A.~Rogers, H.~Groemer and H.~Zassenhaus. They also established a
 joint theory of finite and infinite packings (and coverings) in the
 plane. However, for higher dimensions, Thue's approach does not yield a 
 joint theory of finite and infinite packings.

For the most interesting case of (free) finite sphere packings, L.~Fejes T\'oth
formulated in 1975 his famous {\em sausage conjecture}, claiming that for 
dimensions $\geq 5$ and any(!) number of unit balls, a linear arrangement
of the balls, i.e., all midpoints are on a line and two consecutive balls touch each other, 
minimizes the  volume of their convex
hull. Currently, the sausage conjecture has been confirmed for all
dimensions $\geq 42$.  
The sausage conjecture shows already that finite packings 
have nothing or only  little
to do with classical infinite packings. This motivated the natural question to replace or
to generalize Thue's density by a more general density, which permits a joint
theory of finite and infinite packings for all dimensions and which
also explains strange phenomena of finite  packings as so-called {\em
  sausage catastrophes}.  

We will show in this note, that the {\em parametric density} found in
1993 is the right answer to this question; it 
contains, in particular, Thue's density as a special case. As the name indicates,
the parametric density depends on a parameter $\rho>0$,
and it turns out that for small $\rho$ sausage-like arrangements are
optimal, 
whereas for large $\rho$ densest infinite packing arrangements are optimal.   
The parametric density, in particular, allows to obtain results on infinite packings via (limit) results on finite packings.
Still, these 
results are, with a few exceptions,  weaker than their classical infinite
counterparts, but it is a tempting and interesting task
to improve them and we kindly invite the reader to do so.

There is also an analogous theory of parametric densities for finite
and infinite covering problems of convex bodies, introduced in
\cite{BetkeHenkWills1995}. In this survey, however, we just consider packing problems
and for a thorough treatment of finite and covering problems, including
the parametric densities, we refer to the book of K.~Böröczky, 
Jr.~\cite{Boeroeczky2004}. As a general reference to packing and
covering of convex bodies see, e.g., \cite{FejesToth1999, Gruber2007}, and for sphere
packings, e.g., \cite{Cohn2017, ConwaySloane1999, Zong1999}.

The paper is organized as follows: 
After providing the necessary definitions and notations regarding
infinite and finite packings in the next
section, the parametric densities will be introduced  in Section 3. In
Section 4
we will discuss the small parameter range and in Chapter 5 the large
one. Finite lattice packings with respect to parametric densities are
discussed in the last section.


\section{Notations and preliminaries}
We are working in the $d$-dimensional Euclidean space $\R^d$, equipped
with the standard inner product $\ip{\vx}{\vy}=\vx^\trans\,\vy$
for $\vx,\vy\in\R^d$ and Euclidean norm
$\enorm{\vx}=\sqrt{\ip{\vx}{\vx}}$. $\bad=\{\vx\in\R^d
:\enorm{\vx}\leq 1\}$ is the Euclidean (unit) ball centered at the
origin $\vnull$ of radius $1$; its boundary $\bd\bad$ is called unit
sphere and will be denoted by $\sad$. The set of all convex bodies
$K\subset\R^d$ is denoted by $\Kd$, i.e., $K\in\Kd$, if $K$ is convex,
closed, bounded and $\inte(K)$, the interior of $K$  is non-empty. The
dimension of a set $S$ is the dimension of its affine hull
$\aff(S)$ and it will be denoted by $\dim S$. 
For $K\in\Kd$ let
\begin{equation} 
 \packk=\{C\subset \R^d: \inte(\vx_i+K)\cap
                             \inte(\vx_j+K)=\emptyset,\,
                             \vx_i\ne \vx_j\in C\}
                           \label{eq:packset}  
                           \end{equation} 
be the set of all  {\em packing sets} of $K$.                           
For $C\in\packk$, the arrangement $C+K$ is  called a  {\em packing} of
$K$. In order to define the density of such a packing we denote by $\vol(S)$ the volume, i.e., the $d$-dimensional
Lebesgue measure of a measurable  set $S\subset\R^d$.  For a finite set
$S\subset\R^d$ its cardinality is denoted by $\# S$, and let
$W^d=[-1,1]^d$ be the cube of edge length $2$ centered at the origin.
Then for $K\in\Kd$, $C\in\packk$, 
\begin{equation}
  \delta(K,C)=\limsup_{\lambda\to\infty}
  \frac{\#(C\cap\lambda\,W^d)\,\vol(K)}{\vol(\lambda\,W^d)}
\label{eq:infden}   
\end{equation}
is called  the {\em density  of  the packing   $C+K$ } and
\begin{equation}
  \deninf=\sup\{\delta(K,C):C\in\packk\}
  \label{eq:infdenopt}
\end{equation}
is called the {\em density of a densest packing of $K$}.

Obviously,
for any finite packing set $C$ we have $\delta(K,C)=0$. The idea of the
quantity $\delta(K,C)$ is to measure how much of the space $\R^d$  is
occupied by $C+K$, i.e., we would like to measure
$\vol(C+K)/\vol(\R^d)$,  and we do it mathematically by approximating
$\R^d$ via the sequence $\lambda\,W^d$. In particular, $\delta(K,C)$
may depend on the gauge body (here $W^d$) by which we approximate $\R^d$.  
It was shown
by Groemer \cite{Groemer1963a}, however,  that the definition of
$\deninf$ is independent of this gauge body, and that there exists an {\em optimal packing
set} $C_K\in\packk$ such that 
\begin{equation*}
  \deninf=\delta(K,C_K)=\lim_{\lambda\to\infty} 
  \frac{\#(C_K\cap\lambda\,W^d)\vol(K)}{\vol(\lambda\,W^d)}. 
\end{equation*}

Now we turn to finite (free) packings and to this end we consider for
an integer $n\in\N$,
\begin{equation*}
  \packkn=\{C\in\packk: \# C\}
\end{equation*}
the set of all \emph{packing sets of cardinality $n$}. Here 
we want to find a packing set $C_{K,n}\in\packkn$ minimizing
$\vol(\conv\,C+K)$ among all $C\in\packkn$, where $\conv$ denotes the
convex hull. Hence, in analogy to \eqref{eq:infden},
\eqref{eq:infdenopt} we denote for $K\in\Kd$ and $C\in\packk$ with  $\#
C<\infty$ by
\begin{equation}
  \delta_1(K,C)=\frac{\# C\,\vol(K)}{\vol(\conv\,C+K)}
\label{eq:finden}  
\end{equation}
the {\em density of the finite packing  $C+K$} and 
\begin{equation}
  \delta_1(K,n)=\sup\{\delta_1(K,C): C\in\packkn\}
\label{eq:findenopt}  
\end{equation}
is called the {\em density of a densest $n$-packing of $K$}. The role of the
index $1$ will become clear soon, and it not hard to see that for any
$n$ there exists an \emph{optimal finite packing set} $C_{n,K}$ such that
$\delta_1(K,n)=\delta_1(K,C_K(n))$. 

Of particular interest are here finite packing sets 
$C=\{\vx_1,\dots,\vx_n\}\in\packkn$ with $\dim C=1$, i.e., all points are collinear.
Since we also want to  minimize $\vol(\conv \{\vx_1,\dots,\vx_n\}+K)$
we may assume that for two consecutive  points on
this line,  $\vx_i,\vx_j$, say, the translates $\vx_i+K$ and $\vx_j+K$
touch. Hence,  without loss of generality the points of such a
packing set can be represented as  
\begin{equation*}
   S_n(K,\vu)=\left\{(i-1)\frac{2}{\norm{\vu}{K}}\vu : 1\leq i\leq n\right\},
 \end{equation*}
 where $\vu\in \sad$ is the direction of the line and with $\norm{\vu}{K}$ we
 denote the norm induced by the origin symmetric body
 $\frac{1}{2}(K-K)$, i.e.,
 \begin{equation*}
   \norm{\vu}{K}=\min\left\{\mu\in\R_{\geq 0} : \vu\in\mu\, \frac{1}{2}(K-K)\right\}.
 \end{equation*}
 
 \begin{figure}[htb]
   \includegraphics{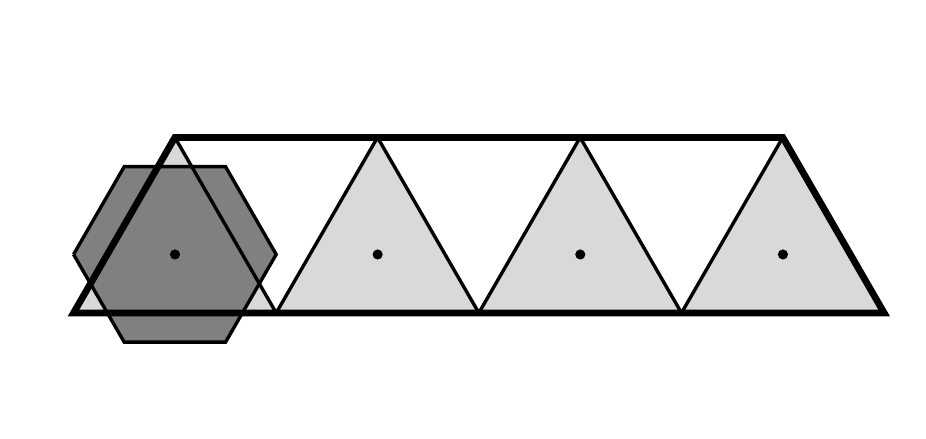}
 \caption{A sausage configuration of a triangle $T$, where
   $\frac{1}{2}(T-T)$  is the darker hexagon.}  
 \end{figure} 

 Packing sets $S_n(K,\vu)$ are called \emph{sausage configurations},
 where the name was coined by L.~Fejes T\'oth \cite{HajnalToth1975a}
 in the special setting $K=\bad$. 
 Obviously, in this case the density of such a sausage
 configuration is independent of the direction $\vu$ and therefore, it
 will be only denoted by $S_n(\bad)$ and it is
 \begin{equation*}
   \vol(\conv(S_n(\bad)+\bad)=2(n-1)\kappa_{d-1}+\kappa_d, 
 \end{equation*}
 where $\kappa_i$ denotes the $i$-dimensional volume of the
 $i$-dimensional unite ball.
 \begin{figure}[htb]
   \begin{center}
     \begin{center}
       \includegraphics{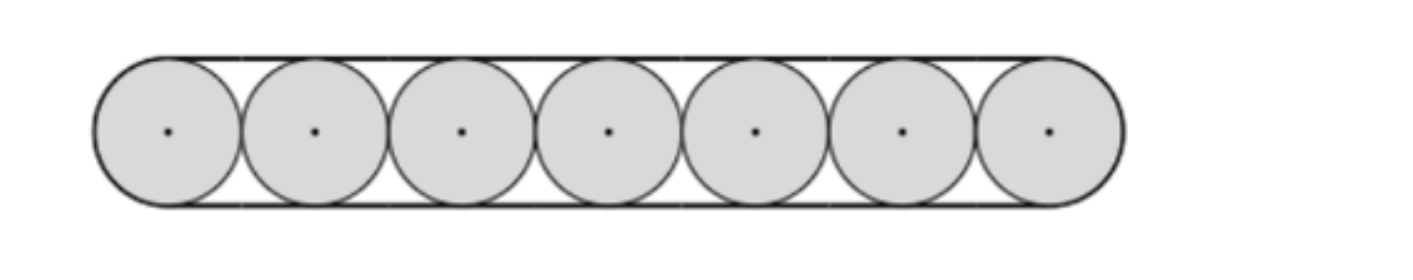}
       \end{center}
     \end{center}
 \caption{A sausage of 7 circles with density $\delta_1(B^2,S_7(B^2))=7\pi/(24+\pi)$.}  
\end{figure} 

The famous sausage conjecture of L.~Fejes T\'oth \cite{HajnalToth1975a}
claims that for any number of balls, a sausage configuration is
 always best possible, provided $d\geq 5$.

 \begin{conjecture}[\sc Sausage conjecture:]  For $d\geq 5$ and $n\in\N$  
 \begin{equation*}
       \delta_1(\bad,n) = \delta(\bad, S_n(\bad)).
     \end{equation*}
 \label{conj:sausage}    
 \end{conjecture}

In the plane  a sausage is never optimal for $n\geq 3$ and
for ``almost all'' $n\in\N$ optimal packing configurations are known
(see \cite{Kenn2011, Schuermann2000, Wegner1986a} and the references within).

In dimension 3 and 4 the situation is more complicated:
In \cite{BetkeGritzmann1984, BetkeGritzmannWills1982}   it was shown that among those finite packings sets
$C$ satisfying 
$\dim C\leq \min\{9,d-1\}$ or $\dim C\leq (7/12)(d-1)$ only sausages are optimal. Hence, in
particular,   for dimensions 3 and 4, no  packings sets of
intermediate dimensions  
are optimal, i.e., optimal packings sets are either $1$-dimensional
(sausages) or $d$-dimensional (clusters).
It is easy to see that for small $n$ sausages are optimal while for
large $n$ clusters are optimal and so the  interesting question is:
when, i.e., for which \emph{''magic'' numbers} $n$ does it happen? 
In dimension 3, results of Wills \cite{Wills1983b, Wills1985c},
Gandini and Wills \cite{GandiniWills1992} and
Scholl\cite{Scholl2000} show that certain
clusters are denser than sausage configurations when
$n=56$ or  $n\geq 58$. In fact, it is conjectured that for $n<56$ and
$n=57$ sausages are optimal.  In dimension 4 it was shown  by
Gandini and Zucco \cite{GandiniZucco1992}, Gandini \cite{Gandini1995b} that a cluster is better than a sausage configuration
for $n\geq 375{,}769$. This large number of spheres 
motivated the  name {\em sausage catastrophe}  given in \cite{Wills1985c}
referring to the abrupt change of the optimal shape of an optimal 
packing set. For a German popular science article about the catastrophe and
the conjecture see \cite{Freistetter2019}.

Obtaining a unified theory for finite and infinite packings covering 
also these  phenomena of sausage conjecture and sausage catastrophe was one 
motivation for the parametric density which we will define in the next
section.

L.~Fejes T\'oth's sausage conjecture was first proved via the parametric
density approach in dimensions $\geq
13{,}387$ by Betke et al.~\cite{BetkeHenkWills1994} which was later
improved to $d\geq 42$ by Betke and Henk  \cite{BetkeHenk1998}. 

The sausage conjecture, in particular, implies that in general 
\begin{equation*}
   \delta(K)< \limsup_{n\to\infty }\delta_1(K,n)
 \end{equation*}
 and, in fact, this is known to be true for all dimensions $d\geq 3$. Thus, large
 optimal finite packing sets do not ``approximate''  optimal  infinite
 packing sets.  However, as we will see next, this will be corrected via the parametric density.    

 \section{The parametric density} 
 The concept of a parametric density was introduced by Betke et al.~in
 \cite{BetkeHenkWills1994}, and the definitions and results presented in
 this section are taken from this paper. 
\begin{definition}[Parametric Density] Let $\rho>0$ and $K\in\Kd$. 
  \begin{enumerate}
    \item  Let $C\in\packk$ with  $C<\infty$, 
  \begin{equation}
    \delta_\rho(K,C)=\frac{\# C\,\vol(K)}{\vol(\conv C +\rho\,K)}
\label{eq:paradens} 
  \end{equation}
  is called the parametric density of $C$ with respect to $K$ and
  the parameter $\rho$.
\item
  \begin{equation*}
    \delta_\rho(K,n)=\sup\{ \delta_\rho(K,C): C\in\packkn\}
  \end{equation*}
   is called the parametric density of a densest $n$-packing of $K$ with respect
   to the parameter $\rho$.  
 \item 
   \begin{equation*}
    \delta_\rho(K)=\limsup_{n\to\infty} \delta_\rho(K,n)
  \end{equation*}
  is called the parametric limit density of $K$ with respect
   to the parameter $\rho$.  
 \end{enumerate}   
\end{definition} 
Apparently, for $\rho=1$,  the definitions in i) and ii) coincide with
the previous given definitions of $\delta_1(K,C)$ \eqref{eq:finden}
and $\delta_1(K,n)$ \eqref{eq:findenopt}, and again it is easy to see
that the $\sup$ in ii) may be replaced by a $\max$.

It is also easy to check, that  $\delta_\rho(K,C)$, $\delta_\rho(K,n)$ and
$\delta_\rho(K)$ are monotonously decreasing and continuous in
$\rho$. Moreover, $\delta_\rho(K,n)$ and $\delta_\rho(K)$ are invariant with
      respect to regular affine transformations of $K$.  By
      calculating the finite parametric density of large clusters of a
      densest infinite packing one gets for all $\rho> 0$ 
      \begin{equation}
        \delta_\rho(K)\geq \delta(K).
        \label{eq:critrho}
      \end{equation}
In order to understand the role of the parameter $\rho$ we briefly
recall some basic facts about mixed volumes, which implicitly appear
in the  denominator of \eqref{eq:paradens}, and for a detailed
account we refer, e.g., to \cite{Gruber2007, Schneider2014}. It is a classical fact
from Convex Geometry, that the  volume of $\conv C+\rho K$ can be
written as a polynomial in $\rho$ of degree $d$, the so-called
\emph{generalized Steiner-polynomial}:
\begin{equation*}
  \vol(\conv C+ \rho K)=\sum_{i=0}^d \binom{d}{i} \rho^{i}\V_i(\conv C;K),
\end{equation*}
where the coefficients $\V_i(\conv C;K)$ are called {\em mixed
volumes}. In particular, we have  $\V_d(\conv C;K)=\vol(K)$,
$\V_0(\conv C;K)=\vol(\conv C)$ and for $K\in\Kd$ it is 
    $\V_i(\conv C;K) = 0$ if and only if  $\dim C < d-i$. 
  So in order to determine $\delta_\rho(K,n)$ we have to minimize
  \begin{equation*}
    \sum_{i=0}^d \binom{d}{i} \rho^{i}\frac{\V_i(\conv C;K)}{\vol(K)} 
  \end{equation*}
 over all $C\in\packkn$. A large parameter $\rho$ gives a strong
 weight on the mixed volumes with a high index $i$. So it seems
 preferable to make mixed volumes with a small index rather large,
 which means that 
 optimal packing set should be of dimensions $d$. Hence,  we can
 expect 
 that for large $\rho$ and large $n$, optimal finite parametric
 densities converge to the density of a densest infinite 
 packing. Therefore, we define
  \begin{definition}[Critical Parameter] For $K\in\Kd$ let
   \begin{equation*}
      \rho_c(K)=\inf\{\rho>0 : \delta_\rho(K)=\delta(K)\}
    \end{equation*}
  \end{definition}  
\begin{figure}[htb]
   \includegraphics[height=3cm]{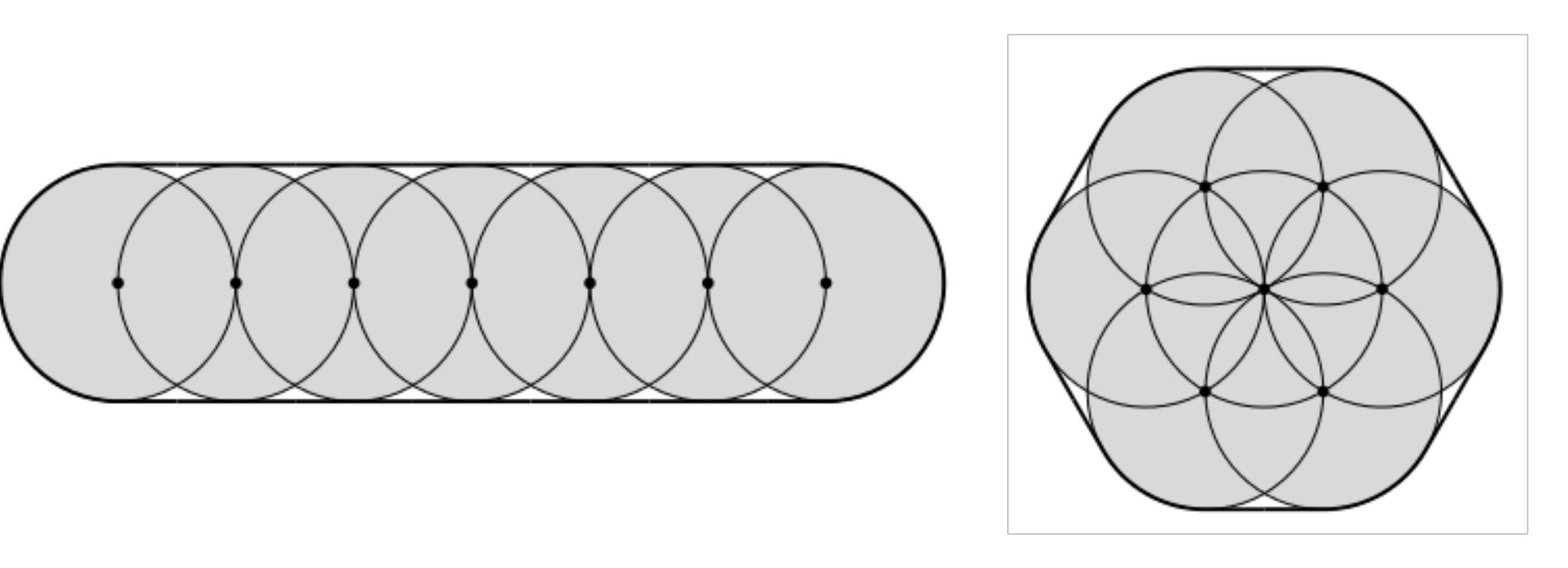}
 \caption{A sausage and a hexagonal packing $C$ of 7 unit circles
   with $\rho=2$. The densities are $\delta_2(B^2,S_7(B^2))=\frac{
     7\pi}{ 2\cdot 24+2^2\cdot\pi}\sim 0.36$ and $\delta_2(B^2,C)=\frac{7\pi}{ 6\sqrt{3}+2\cdot  12+2^2\cdot\pi}\sim
 0.47$.}  
 \end{figure} 

 On the other hand, a small $\rho$ gives a strong weight to the mixed
 volumes with a small index, and so low-dimensional packing sets $C$
 are better, with the extreme case $\dim C=1$. Now for such a sausage
 configuration $S_n(\vu,K)$ we have
 \begin{equation}
   \vol(\conv S_n(\vu,K)+ \rho K) =
   2(n-1)\frac{\vol_{n-1}(K|\vu^\perp)}{\norm{\vu}{K}}
   \rho^{d-1}+\vol(K)\rho^d, 
\label{eq:sausage} 
 \end{equation}
 where $\vol_{n-1}(K|\vu^\perp)$ is the $(n-1)$-dimensional volume of
 the orthogonal porjection of $K$ on the hyperplane orthogonal to
 $\vu$. Hence, in order to have a sausage configuration of maximal
 density let $\vu_k\in S^{n-1}$ be such that
 \begin{equation*}
   \frac{\vol_{n-1}(K|\vu^\perp_K)}{\norm{\vu_K}{K}}=\min_{\vu\in S^{d-1}}\frac{\vol_{n-1}(K|\vu^\perp)}{\norm{\vu}{K}}.
 \end{equation*}
 With this notation let
 \begin{equation}
   \delta_\rho^s(K)= \lim_{n\to\infty} \delta_\rho(K,S_n(\vu_K,K)) =
   \rho^{1-d}\vol(K)\frac{\norm{\vu_K}{K}}{2\,\vol_{n-1}(K|\vu^\perp_K)}
   \label{eq:limitsau}
 \end{equation}
 be the \emph{parametric limit density} of an optimal sausage configuration of
 $K$ with respect to the parameter $\rho$. This density may be
 regarded as the $1$-dimensional counterpart to $\delta(K)$. In
 particular, for $K=\bad$ we have
 \begin{equation}
   \delta_\rho^s(\bad)=\rho^{1-d}\frac{\kappa_{d}}{2 \kappa_{d-1}}. 
\label{eq:ballsaulimit}
 \end{equation}
In analogy to the critical parameter we define
 \begin{definition}[Sausage Parameter] For $K\in\Kd$ let
   \begin{equation*}
      \rho_s(K)=\sup\{\rho>0 : \delta_\rho(K)=\delta_\rho^s(K)\}.
    \end{equation*}
  \end{definition}  

  \begin{figure}[hbt] 
  \includegraphics[height=3cm]{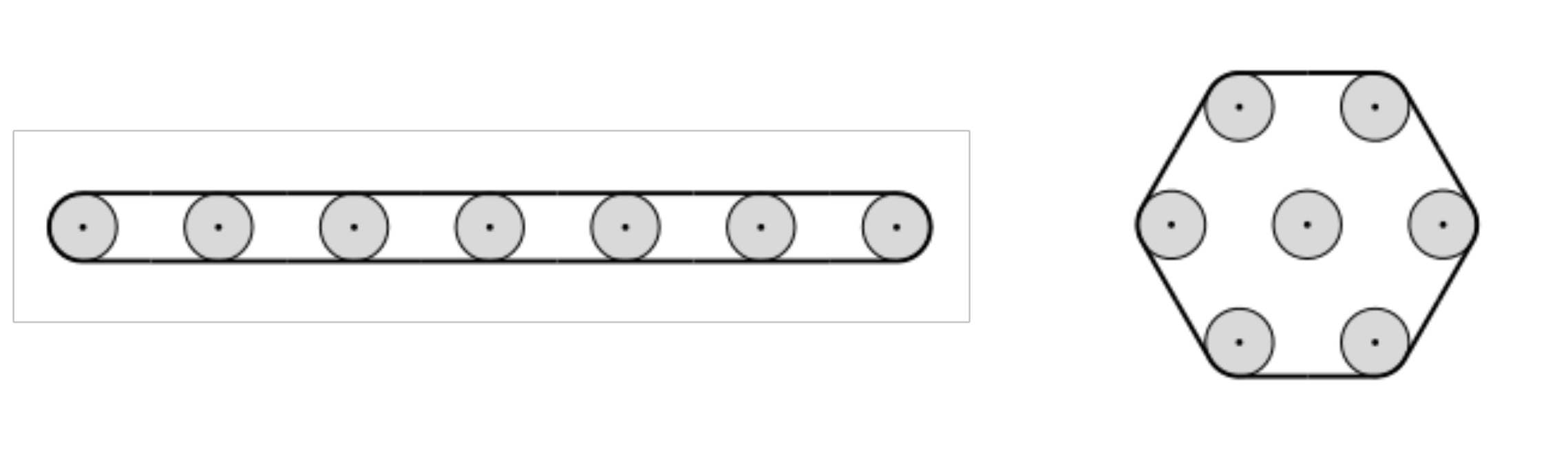}
  \caption{A sausage and a hexagonal packing $C$ of 7 unit circles
    with $\rho=1/2$. The densities are
    $$ 
    \delta_{\frac{1}{2}}(B^2,S_7(B^2))=\frac{7\pi}{(1/2)\cdot
            24+(1/2)^2\cdot\pi}\sim 1.72,$$
    $$     \delta_{\frac{1}{2}}(B^2,C)=\frac{7\pi}{6\sqrt{3}+(1/2)\cdot 12+(1/2)^2\cdot\pi}\sim 1.28.$$
     }  
  \end{figure} 

 These two parameters devide the range of all parameters into three
 relevant areas and provide us with an simple bound on $\delta(K)$
 as  the next theorem shows.
 \begin{theorem}Let $K\in\Kd$ and assume that
   $0<\rho_s(K),\rho_c(K)<\infty$. Then
   \begin{enumerate}
   \item $\rho_s(K)\leq\rho_c(K)$.   
   \item $\delta_\rho(K)=\delta_\rho^s(K)$ for $\rho\in(0, \rho_s(K)]$.
   \item $\delta_\rho(K)=\delta(K)$ for $\rho\in[\rho_c(K),\infty)$.
   \item $\delta_{\rho_c(K)}^s(K)\leq\delta(K)\leq\delta_{\rho_s(K)}^s(K)$.   
   \end{enumerate}
 \label{theo:basic}   
 \end{theorem}   
 In other words, below the sausage parameter infinite sausages, i.e.,
 1-dimensional packings,  are
 optimal, above the critical parameter densest infinite packing
 arrangements are optimal, and what happens in between is (in general)
 rather  unknown. Moreover, the parametric density of infinite
 sausages with respect to  critical and sausage parameter yields lower
 and upper bounds on $\delta(K)$ (cf.~\eqref{eq:limitsau})
 \begin{equation}
        \rho_c(K)^{1-d}\frac{\norm{\vu_K}{K}\vol(K)}{2 \vol_{n-1}(K\vert \vu_K^\perp)}\leq
        \delta(K)\leq \rho_s(K)^{1-d}\frac{\norm{\vu_K}{K}\vol(K)}{2
          \vol_{n-1}(K\vert \vu_K^\perp)}.
 \label{eq:denfininfgen}     
    \end{equation} 
 In particular, for the unit ball
 $\bad$, Theorem \ref{theo:basic} iv) gives
 (cf.~\eqref{eq:ballsaulimit})
 \begin{equation}
     \frac{\kappa_d}{ 2\kappa_{d-1}} (\rho_c(\bad))^{1-d}\leq
      \delta(\bad)\leq
      \frac{\kappa_d}{ 2\kappa_{d-1}} (\rho_s(\bad))^{1-d}.
 \label{eq:ballbounds}     
\end{equation}
Here useful estimates are (see, e.g., \cite{Gritzmann1985c}, \cite{BetkeGritzmannWills1982})
\begin{equation}
  \frac{1}{d}<\frac{\norm{\vu_K}{K}\vol(K)}{2 \vol_{n-1}(K\vert
    \vu_K^\perp)}\leq 1\text{ and }
  \sqrt{\frac{2\pi}{d+1}}<\frac{\kappa_d}{\kappa_{d-1}}<\sqrt{\frac{2\pi}{d}}.
 \label{eq:boundballvolume} 
\end{equation}

In general, we have $\rho_s(K)<\rho_c(K)$ for $d\geq 3$ (cf.~\cite[Theorem 10.7.1]{Boeroeczky2004}), but in
case of the ball it is tempting to conjecture that
$\rho_c(\bad)=\rho_s(\bad)$. In fact, for the ball it is even
conjectured in \cite{BetkeHenkWills1994}
\begin{conjecture}[\sc Strong Sausage Conjecture] For $n\in\N$ and
  $\rho>0$ it holds
  \begin{equation*}
     \delta_\rho(\bad,n)=\delta_\rho(\bad,S_n(\bad)) \text{ or } \delta_\rho(\bad,n)<\delta(\bad).
   \end{equation*}
In particular, we have $\rho_c(\bad)=\rho_s(\bad)$.   
\label{conj:ssc}  
\end{conjecture}
Let us briefly point out that this conjecture also covers
(essentially) L.~Fejes
T\'oth's sausage conjecture as for $\rho=1$ and $d\geq 5$ a (maybe
large) sausage has a larger density than $\delta(\bad)$.   
This conjecture  would also imply the equivalence of determining $\delta(\bad)$,
$\rho_c(\bad)$ and $\rho_s(\bad)$.  
 In the next sections we will briefly present what is known about these
 two parameters.

 \section{The planar case}
 In the planar case and for centrally symmetric convex domains
 we have a rather complete picture.  
 Based on classical results of Rogers \cite{Rogers1951a, Rogers1960a}
 and  Oler \cite{Oler1961a}  it was
 shown in \cite{BetkeHenkWills1994} that 
 \begin{theorem} Let $K\in\mathcal{K}^2$, $K=-K$, and $n\in \N$. 
   \begin{enumerate}
\item
 \begin{equation*}               
                \frac{3}{4}\leq \rho_s(K)=\frac{\delta_1^s(K)}
  {\delta(K)}=\rho_c(K)\leq 1,
\end{equation*} 
where equality on the left  is attained only  for an affinely regular
non-degenerate hexagon and on the right only for a parallelogram.      
\item $\delta_\rho(K,n)=\delta_\rho(K,S_n(K))$ for $0<\rho\leq \rho_s(K)$.
\item $\delta_\rho(K,n)\leq \delta(K)\left(\frac{n}
       { n-1+\delta(K)\rho^2}\right)$,\quad $\rho_s(K)\leq\rho<\infty$.
   \end{enumerate}
\label{thm:planar}    
 \end{theorem}

For instance, for the circle $B^2$ we have $\delta_1^s(B^2)=\pi/4$
(cf.~\eqref{eq:ballsaulimit}) and by a well-known and classical result of
Thue \cite{Thue1892} we also know $\delta(B^2)=\pi/(2\sqrt{3})$. Hence
\begin{equation*}
  \rho_s(B^2)=\frac{\sqrt{3}}{2}.
\end{equation*} 

In particular, Theorem \ref{thm:planar} ii) shows that for $\rho<\rho_s(K)$ and
any $n\in\N$ only sausage configurations are optimal. Such a strong
statement is not true for arbitrary convex domains in the plane as it
was  shown by Böröczky, Jr. and Schnell in \cite{BoeroeczkySchnell1998}, where they also
proved that  Theorem \ref{thm:planar} i) holds true for any convex body in the
plane. 

For $\rho \geq \rho_s(K)$ various optimal configurations might be
possible as our running example below shows. 

 \begin{figure}[htb]
   \includegraphics[height=3cm]{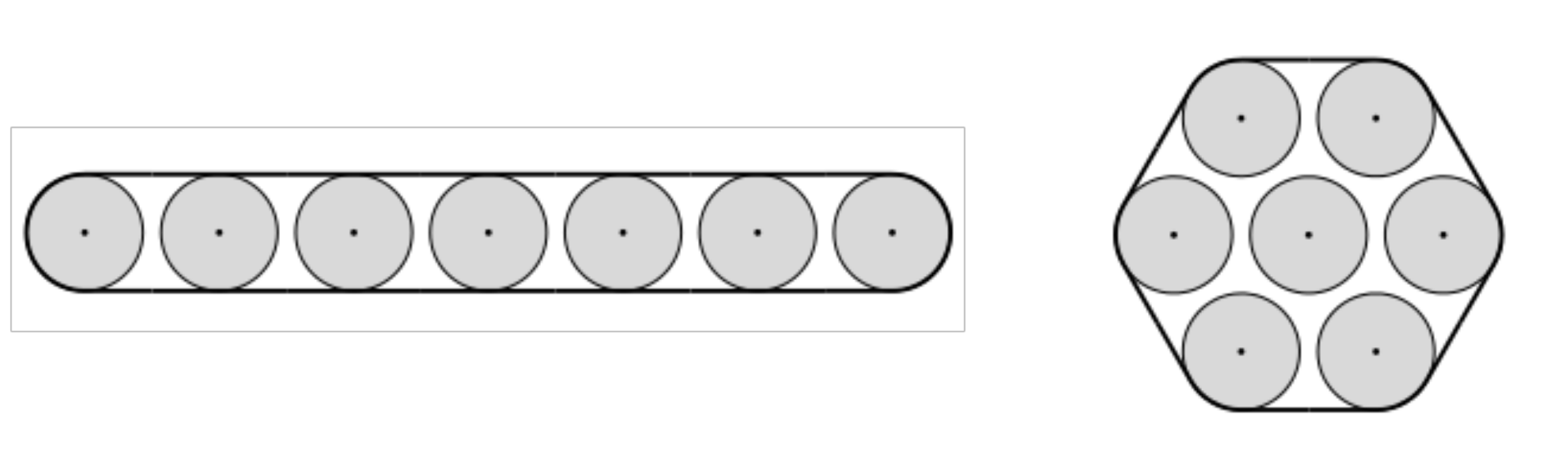}
 \caption{A sausage and a hexagonal packing $C$ of 7 unit circles
   with $\rho=\rho_s(B^2)=\sqrt{3}/2$. Here we have
   $\delta_{\sqrt{3}/2}(B^2,S_7(B^2))=\delta_{\sqrt{3}/2}(B^2,C)=7\,\pi/(6\sqrt{3}+12\sqrt{3}/2+\pi\,(\sqrt{3}/2)^2)\sim
   0.9503$.}  
 \end{figure} 
 
\section{Small parameters and sausages}
First we focus on the most prominent convex body, the ball. Here the
main result verifies the sausage conjecture in high dimensions in a very strong form,
i.e., not
only for $\rho=1$. 
\begin{theorem}[\cite{BetkeHenkWills1995a}] For every $\rho<\sqrt{2}$
  there exists a dimension $d_\rho$ such that for all $n\in\N$ and
  $d\geq d_\rho$
  \begin{equation*}
    \delta_\rho(\bad,n)=\delta(\bad,S_n(\bad)).
  \end{equation*}
\label{thm:smallball}   
\end{theorem}
Even for parameters $\geq 1$ sausages are optimal packing
configurations for any $n\in\N$, provided the dimension is large
enough. Of course, this result implies
\begin{corollary}  $\liminf_{d\to\infty}\rho_s(\bad)=\sqrt{2}$.
\end{corollary} 

The proof of Theorem \ref{thm:smallball} is based on a local approach
by comparing the volume of $\conv(C+\rho\bad)$ contained in a Dirichlet-Voronoi cell of an arbitrary packing set $C\in\packbn$ to the
corresponding volumes of a sausage  $S_n(\bad)+\rho\bad$. The bound of
$\sqrt{2}$ corresponds  to the minimum distance between a vertex of
a Dirichlet-Voronoi cell at $c$, say, and the center $c$. It was shown
by Rogers \cite[Chapter 7]{Rogers1964a} 
that the distance between an $i$-dimensional face of such a cell and its center is at
least $\sqrt{2(n-i)/(n-i+1)}$.  In fact, in \cite{BetkeHenkWills1994}, Theorem
\ref{thm:smallball} was  proved for $\rho<2/\sqrt{3}$ since
only the $(n-2)$-dimensional faces of the cell were ``used.'' It is
easy to see that not all vertices of a Dirichlet-Voronoi cell of a
packing set $C\in\packbn$ can be as close as  $\sqrt{2}$, but it seems
to be hard to take advantage of this fact.

Another consequence of Theorem \ref{thm:smallball} is the following upper bound on
$\delta(\bad)$ (cf.~\eqref{eq:ballbounds} and
\eqref{eq:boundballvolume}).
\begin{corollary} For $\epsilon>0$ there exists a dimension
  $d_\epsilon$ such that for $d\geq d_\epsilon$
  \begin{equation*}
    \delta(\bad)\leq \sqrt{\frac{\pi}{d}}(\sqrt{2}-\epsilon)^{1-d}.
  \end{equation*} 
\end{corollary}
This   bound   is 
asymptotically  of  the  same  order   as  the  classical  bounds  of  
Blichfeldt  \cite{Blichfeldt1929a}  and  Rogers \cite[Chapter 7]{Rogers1964a}.
Though  this  is   much  weaker  than  the  best   known  upper  bound  for
$\delta(B^d)$  (cf.~\cite{KabatjanskiuiLevensteuin1978}),  it shows that finite
parameterized packings are also a tool to study infinite packings.

In view of L.~Fejes T\'oth's sausage conjecture the dimension $d_1$
of the theorem above is of particular interest. 
\begin{theorem}[\cite{BetkeHenk1998}] $d_1\leq 42$, i.e., the Sausage Conjecture \ref{conj:sausage} is true for all dimension $\geq 42$.   
\end{theorem}
The reason for $42$ is given here \cite{Adams1979}. The sausage conjecture has also
been verified with respect to certain restriction on the packings
sets, e.g., among those which are lower-dimensional \cite{BetkeGritzmann1984, BetkeGritzmannWills1982}, or close to
sausage-like  arrangements \cite{KleinschmidtPachnerWills1984}, of whose inradius is rather large
\cite{BoeroeczkyHenk1995}. For detailed information we refer to \cite[Section 8.3]{Boeroeczky2004}.

Now regarding the sausage parameter of arbitrary convex bodies the
best bound is due to K.~Böröczky, Jr.; he showed
\begin{theorem}[\protect{\cite[Theorem 10.1.1]{Boeroeczky2004}}] Let $K\in\Kd$ and let $\rho<1/(32\,d)$.
  Then for all $n\in\N$ 
  \begin{equation*}
    \delta_\rho(K,n)=\delta(\bad,S_n(K)).
  \end{equation*}
\label{thm:smallbody}   
\end{theorem}
A weaker upper bound of $1/(32\,d^2)$ was frist proved in
\cite{BetkeHenkWills1995a}. By the theorem above we get 
\begin{corollary} Let $K\in\Kd$. Then $\rho_s(K)\geq\frac{1}{32d}$. 
\end{corollary}
It is tempting to conjecture that $\rho_s(K)$ is bounded from below by
an absolute constant.  

\section{Large parameters and densest infinite packings}

Here the main result is captured by the next  theorem. 
\begin{theorem} Let $K\in\Kd$ and let $\ov\rho\in\R_{>0}$ such that
  $K-K\subseteq\,\ov\rho K$. Then for each $n\in\N$
\begin{equation*} 
  \delta_{\ov\rho}(K,n)\leq\delta(K).
\end{equation*}
\label{thm:largebody} 
\end{theorem}
 Hence, for such a $\ov\rho$ we have $\rho_c(K)\leq \ov\rho$
 (cf.~\eqref{eq:critrho}). In order to estimate $\ov\rho$ we may assume by the
 translation invariance of the densities that the
 centroid of $K$ is at the origin. Then, it is known that $K-K\subseteq
 (d+1)K$ and if $K$ is origin symmetric then, of course, $K-K=2\,K$.
 So we have
 \begin{corollary} Let $K\in\Kd$. Then
   \begin{equation*}
     \rho_c(K)\leq\begin{cases}2 &: K=-K\\
                              d+1 &: \text{otherwise}\end{cases}. 
   \end{equation*}
 \end{corollary}
The general bound of $d+1$ was slightly improved to $\rho_c(K)\leq d$ in
\cite[Lemma 10.5.2]{Boeroeczky2004}, but also here it might be true that
$\rho_c(K)$ is bounded from above by an absolute constant. 

The proof of Theorem \ref{thm:largebody} is based on an average
argument: We assume that there exists a $C\in\packkn$ with
$\delta_{\ov\rho}(K,n)>\delta(K)$. For each $\vx\in[0,\gamma]^n$, for some
large $\gamma$,  a finite
packing set $C_\vx$ contained in a large cube $\lambda\,W^d$ is
constructed consisting of suitable translates of $C$ and of points of a densest infinite packing of $K$. For this
superposition the property $K-K\subset\ov\rho\,K$ is used.
In order to determine the cardinality of this new packing set $C_\vx$
it is averaged over all $\vx\in  [0,\gamma]^n$. This yields the
existence of a packing set $C_{\ov\vx}$ of $K$ such that $\#
C_{\ov\vx}\vol(K)\ \vol(\lambda\,W^d)>\delta(K)$, which contradicts the
definition of  $\delta(K)$.  

Regarding infinite packings, the corollary above 
together with \eqref{eq:denfininfgen} with \eqref{eq:ballbounds} gives for $K\in\Kd$, $K=-K$, the
bound 
\begin{equation*}
  \delta(K)> \frac{1}{d}\,2^{1-d}.
\end{equation*}
This is up to a factor of $1/d$ of the same order
than the best known lower bounds  on $\delta(\bad)$ (cf., e.g., \cite[Theorem
2.2]{Rogers1964a}).

In order to present results regarding the shape of optimal finite
packings, we denote  for $\rho\in\R_{>0}$ and   $n\in
  \N$ by $C_{\rho,n,K}$ an \emph{optimal $n$-packing set of $K$ with respect to the
  parameter} $\rho$. 
  Schnell and Wills \cite{SchnellWills2000} proved that optimal sphere packing
  sets are never two-dimensional. More precisely, 
  \begin{theorem} For any $d\geq 3$, $\rho>0$ and $n\geq 4$ we have $\dim
    C_{\rho,n,\bad}\ne 2$.
  \label{thm:flatball}  
\end{theorem}
A direct consequence is that for $d=3$ and all $\rho\in(\rho_c,1)$ 
 sausage catastrophes occur for $n_\rho\geq 56$ balls.
Exact values of these ''magic'' numbers $n_\rho$ are not known, but they grow when $\rho$
approches $\rho_c$.

For large  $\rho$ and arbitrary convex bodies 
it was shown by Böröczky, Jr. and Schnell \cite{BoeroeczkySchnell1998} that  $\conv C_{\rho,n,K}$
resembles the shape of $K$; more precisely: 
\begin{theorem} Let $K\in\Kd$ and 
  \begin{enumerate}
  \item Let $\rho>\rho_c(K)$. Then
    \begin{equation*}
      \lim_{n\to \infty} \mathrm{r}_K(\conv C_{\rho,n,K})=\infty,
    \end{equation*}
    where $\mathrm{r}_K(\conv C_{\rho,n,K})$ is the (relative) inradius with
    respect to $K$. i.e.,  $\mathrm{r}_K(\conv C_{\rho,n,K})=\max\{r :
                  \vt+rK\subset\conv C_{\rho,n,K}\text{ for some
                  }\vt\in\R^d\}$.
                \item Let $\rho > d+1$. There exists a constant
                  $\mu(\rho,d)$ such that for large $n$  (and after
                  possible translation) 
                  \begin{equation*}
                   \frac{1}{\mu(\rho,d)} K\subset
                   \sqrt[d]{\frac{\delta(K)}{n}}\,C_{n,\rho,K} \subset
                   \mu(\rho,d) K. 
                 \end{equation*}
                 Moreover, $\lim_{\rho\to \infty}\mu(\rho,d)=1$, and 
                 if $K=-K$ the bound on $\rho$ can be lowered to $2$.
               \end{enumerate}
   \label{thm:shape}            
\end{theorem}

\section{(Finite) lattice packings}
In this section we will restrict the finite and infinite packings to
lattice packings. Lattice packings   of convex bodies have a long and famous history
for which we refer to \cite{Gruber2007}. Here we just briefly recall a few facts
which are relevant for finite parameterized lattice packings.

We will
understand by a lattice $\Lambda\subset\R^d$ a regular linear image of
the integral  lattice $\Z^d$, i.e., there exists a regular matrix
$B\in\R^{d\times d}$ such that $\Lambda=B\Z^d$. The determinant of the
lattice, denoted by $\det\Lambda$,  is the volume of the
parallelepiped spanned by the columns of $B$, i.e., $\det\Lambda=|\det B|$.

In analogy to \eqref{eq:packset}
the set of all \emph{packing lattices} of a convex body is denoted by
$\packk^*$, i.e.,
\begin{equation*}
    \packk^*=\{\Lambda\in\packk : \Lambda \text{ lattice}\}.
  \end{equation*}
For $\Lambda\in\packk^*$ the density $\delta(K,\Lambda)$ (cf.~\eqref{eq:infden})
can easily be calculated as $\delta(K,\Lambda)=\vol(K)/\det\Lambda$.
In the case of (infinite) lattice packings  we are interested in the
determination of 
\begin{equation*}
              \delta^*(K)=\sup\{\delta(K,\Lambda): \Lambda\in\packk^*\}.
\end{equation*}
Confirming a conjecture of Minkowski, Hlawka \cite{Hlawka1944} proved a lower
bound on   $\delta^*(K)$ of order to $2^{-d}$, which was subsequently
(slightly) improved by various authors. The current record is still due to
W. Schmidt \cite{Schmidt1963b} with
\begin{equation*}
    \delta^*(K) \geq c\,d\,2^{-d},
\end{equation*}
where $c$ is an absolute constant. For the ball $\bad$ there are even better
 bounds available and here we refer to  the survey of H. Cohn \cite{Cohn2017a}
 and the references within. 

In order to deal with finite lattice packings we restrict  now all given definitions
in Section 3 to lattices, i.e., we set for $K\in\Kd$,
$\rho\in\R_{>0}$ 
\begin{equation*}
   \begin{split}
     \packkn^* & =\{C\in\packkn: \text{there exits a } \Lambda\in\packk^* \text{ with }C\subset\Lambda\} \\   
    \delta_\rho^*(K,n) & =\sup\{\delta_\rho(K,C) : C\in\packkn^*\}\\
    \delta_\rho^*(K) & = \limsup_{n\to\infty}   \delta_\rho^*(K,n)\\ 
    \rho_c^*(K) &= \inf\{\rho\in\R_{>0} : \delta_\rho^*(K)
    =\delta^*(K)\}\\ 
    \rho_s^*(K)&=\sup\{\rho\in\R_{>0} : \delta_\rho^*(K)=\delta_\rho^s(K)\}
   \end{split}
\end{equation*}
and we call all these quantities as in the general case except that we
add the word ``lattice''.  

Then Theorem \ref{theo:basic}  holds also true for all these lattice
quantities, and we also have (cf.~\eqref{eq:critrho})
\begin{equation}
  \delta_\rho^*(K)\geq \delta^*(K).
\label{eq:limlat}   
\end{equation}
Since $\packkn^*\subset\packkn$ it is 
$\rho_s^*(K)\geq \rho_s(K)$, and  so all bounds presented in Section
5 for $\rho_s(K)$ are valid for $\rho_s^*(K)$ as well.

Regarding the critical lattice parameter $\rho_c^*(K)$ the situation is different; in particular, the
superposition argument leading to Theorem \ref{thm:largebody} and to the lower bounds
on $\rho_c(K)$ does not work in the lattice case since the
superposition may destroy the lattice property. 

Here the idea is to use a lattice refinement argument which goes back to
Rogers \cite{Rogers1950}.  It implies  that for $K=-K$ and
$C\subset\packkn^*$ one can always find a lattice 
$\Lambda\subset\packk^*$ containing $C$ such that all points in
space are close to some lattice points; more precisely,  at most at
distance $3$ measured with respect to the norm induced by $K$. Hence, the
circumradius of the Dirichlet-Voronoi cells of this refined lattice $\Lambda$
is at most $3$. Together
with some improvements for the case $K=\bad$,
one gets, 
roughly speaking,  the
following theorem:
\begin{theorem}[\cite{Henk1995}] Let $K\in\Kd$ and $n\in\N$.   Then
\begin{equation*} 
        \delta_\rho^*(K,n)\leq \delta^*(K)
\hbox{ for }
 \rho\geq       \left\{\begin{array}{l@{\,:\,}l}
        \sqrt{21}/2 & K=\bad , \\ 3 & K=-K,
                         \\(3/2)(d+1)& \hbox{otherwise }.
        \end{array}\right.
\end{equation*} 
\end{theorem}
In combination with \eqref{eq:limlat} we obtain
\begin{corollary} 
  Let $K\in\Kd$. Then
\begin{equation*}   
 \rho_c^*(K)\leq      \left\{\begin{array}{l@{\,:\,}l}
        \sqrt{21}/2 & K=\bad, \\ 3 &  K=-K, \\(3/2)(d+1)& \hbox{otherwise }.
                             \end{array}\right.
\end{equation*}                           
\end{corollary}
For a detailed account to Rogers lattice refinement method as well as improvements for $l_p$-ball packings we refer to
\cite{Henk2018} and the references within.

Regarding the structure of optimal finite lattice packings we first
point out that Theorem \ref{thm:flatball} has also been proven by
Schnell and Wills
in the lattice case, i.e.,
\begin{theorem}[\cite{SchnellWills2000}] For any $d\geq 3$, $\rho>0$ and $n\geq 4$ we have $\dim
    C_{\rho,n,\bad}^*\ne 2$.
  \end{theorem}
Regarding optimal packing sets for $d=3$ and $n$ small we refer to
\cite{SchollSchuermannWills2003}. 
  For large $\rho$ and $n$, the asymptotic shape of optimal finite
lattice packings of convex bodies is closely
related to the problem to understand crystal growth and to the so
called \emph{Wulff shape}. We are not going to enter this subject, instead we
refer to \cite[Section 10.11]{Boeroeczky2004}. However, as a strong counterpart to Theorem \ref{thm:shape} i), we
mention a result of Arhelger et al. on origin symmetric convex boides showing that for $\rho$ beyond the
critical lattice parameter the asymptotic shape is cluster like.
\begin{theorem}[\cite{ArhelgerBetkeBoeroeczky2001}] Let $K\in\Kd$, $K=-K$, $\rho>\rho_c^*(K)$. Then the
  ratio of circumradius of $\conv C_{\rho,n}^*$ to the inradius of
  $\conv C_{\rho,n}^*$ is bounded (independent of $n$).
\end{theorem}
In other words, the shape of $\conv C_{\rho,n}^*$ is not far from being a ball.
In the same paper, it was also shown  that for the
3-dimensional ball 
\begin{theorem}[\cite{ArhelgerBetkeBoeroeczky2001}] 
               $\rho_s^*(B^3)=\rho_c^*(B^3)$,
\end{theorem} 
indicating that also a lattice analogue of the Strong Sausage Conjecture \ref{conj:ssc}
could be true (as well). Further evidence is contained in the paper \cite{BoeroeczkyWills1997}.

\medskip
\noindent 
{\it Acknowledgement.} The authors thank the referee for many valuable
comments.




\end{document}